\documentclass[a4paper,12pt]{article}
\usepackage{a4wide}
\usepackage{amsmath}
\usepackage{amssymb}
\usepackage{amsthm}
\usepackage{latexsym}
\usepackage{graphicx}
\usepackage[english]{babel}
\usepackage{makeidx}

\newtheorem{theorem}{Theorem}[section]
\newtheorem{lemma}[theorem]{Lemma}
\newtheorem{corollary}[theorem]{Corollary}

\newtheorem{conjecture}[theorem]{Conjecture}

\theoremstyle{definition}
\newtheorem{definition}[theorem]{Definition}

\theoremstyle{remark}

\begin{document}
\title{The Cones Associated to Some Transvesal Polymatroids}
\author{Alin \c{S}tefan}
\date{}
\maketitle
\begin{abstract}
In this paper we describe the facets cone associated to
transversal polymatroid presented by
${{\bf{\mathcal{A}}}}=\{\{1,2\},\{2,3\},\ldots,\{n-1,n\},\{n,1\}\}.$
Using the Danilov-Stanley theorem to characterize the canonicale
module, we deduce that the base ring associated to this
polymatroid is Gorenstein ring. Also, starting from this
polymatroid we describe the transversal polymatroids with
Gorenstein base ring in dimension 3 and with the help $\it
Normaliz$ in dimension 4.

\vspace{5 pt} \noindent \textbf{Keywords:} Base ring, Transversal
polymatroids, Danilov-Stanley theorem, Equations of the cone.

\vspace{5 pt} \noindent \textbf{2000 Mathematics Subject
Classification:} 05C50,13A30,13H10,13D40.
\end{abstract}

\section{Preliminaries on polyhedral geometry}
An affine space generated by $A\subset {\bf{R^{n}}}$ is a
translation of a linear subspace of ${\bf{R^{n}}}.$ If $0\neq a\in
{\bf{R^{n}}},$ then $H_{a}$ will denote the hyperplane of
${\bf{R^{n}}}$ through the origin with normal vector a, that
is,\[H_{a}=\{x\in {\bf{R^{n}}}\ | \ <x,a>=0\},\] where $<,>$ is
the usual inner product in ${\bf{R^{n}}}.$ The two closed
halfspaces bounded by $H_{a}$ are: \[H^{+}_{a}=\{x\in
{\bf{R^{n}}}\ | \ <x,a> \geq 0\} \ and \  H^{-}_{a}=\{x\in
{\bf{R^{n}}}\ | \ <x,a> \leq  0\}.\]

Recall that a $polyhedral\ cone\ Q \subset {\bf{R^{n}}}$ is the
intersection of a finite number of closed subspaces of the form
$H^{+}_{a}.$ If $A=\{\gamma_{1},\ldots, \ \gamma_{r}\}$ is a
finite set of points in ${\bf{R^{n}}}$ the $cone$ generated by
$A$, denoted by ${\bf{R_{+}}}{A},$  is defined as
\[{\bf{R_{+}}}{A}=\{\sum_{i=1}^{r}a_{i}\gamma_{i}\ | \ a_{i}\in
{\bf R_{+}},\ with \ 1\leq i \leq n\}.\]

An important fact is that $Q$ is a polyhedral cone in
${\bf{R^{n}}}$ if and only if there exists a finite set $A\subset
{\bf{R^{n}}}$ such that $Q={\bf{R_{+}}}{A},$ see ({\cite
{W},theorem 4.1.1.}).

\begin{definition}
A proper face of a polyhedral cone is a subset $F\subset \ Q$ such
that there is a supporting hyperplane $H_{a}$ satisfying:

$1)$ $F=Q\cap H_{a}\neq \emptyset.$

$2)$ $Q\nsubseteq H_{a}$ and $Q\subset H^{+}_{a}.$
\end{definition}

\begin{definition}
A proper face $F$ of a polyhedral cone $Q\subset \ {\bf{R^{n}}}$
is called a $facet$ of $Q$ if $dim(F)=dim(Q)-1.$
\end{definition}

\section{Polymatroids}
Let $K$ be a infinite field, $n$ and $m$ be positive integers,
$[n]=\{1, 2, \ldots , n\}$. A nonempty finite set $B$ of ${\bf
N}^{n}$ is the base set of a discrete polymatroid
${\bf{\mathcal{P}}}$ if for every $u=(u_{1}, u_{2}, \ldots,
u_{n})$, $v=(v_{1}, v_{2}, \ldots, v_{n})$ $\in B$ one has $u_{1}
+ u_{2} + \ldots + u_{n}=v_{1} + v_{2} + \ldots + v_{n}$ and for
all $i$ such that $u_{i}>v_{i}$ there exists $j$ such that
$u_{j}<v_{j}$ and $u+e_{j}-e_{i} \in B $, where $e_{k}$ denotes
the $k^{th}$ vector of the standard basis of ${\bf N}^{n}$. The
notion of discrete polymatroid is a generalization of the
classical notion of matroid, see \cite{E} \cite{HH} \cite{O}
\cite{W1}. Associated with the base $B$ of a discret polymatroid
${\bf{\mathcal{P}}}$ one has a $K-$algebra $K[ B ]$, called the
base ring of ${\bf{\mathcal{P}}}$, defined to be the
$K-$subalgebra of the polynomial ring in $n$ indeterminates
$K[x_{1}, x_{2}, \ldots , x_{n}]$ generated by the monomials
$x^{u}$ with $u \in B$. From \cite{HH} the algebra $K[ B ]$ is
known to be normal and hence Cohen-Macaulay.

If $A_{i}$ are some non-empty subsets of $[n]$ for $1\leq i\leq
m$, ${\bf{\mathcal{A}}}=\{A_{1},A_{2},\ldots,A_{m}\}$, then the
set of the vectors $\sum_{k=1}^{m} e_{i_{k}}$ with $i_{k} \in
A_{k},$ is the base of a polymatroid, called transversal
polymatroid presented by ${\bf{\mathcal{A}}}.$ The base ring of a
transversal polymatroid presented by ${\bf{\mathcal{A}}}$ denoted
by $K[{\bf{\mathcal{A}}}]$ is the ring :
\[K[{\bf{\mathcal{A}}}]:=K[x_{i_{1}}x_{i_{2}}\ldots
x_{i_{m}}:i_{j}\in A_{j},1\leq j\leq m].\]

\section{Some Linear Algebra}
Let $n\in{\bf N}$ a integer number, $n\geq 3$ and the following
set with $2n-3$ points with positive integer coordinates :
\[ \{ \ R_{0,1}(2,1,1,\ldots,1,1,0), R_{0,2}(2,1,1,\ldots,1,0,1),\ldots,
R_{0,n-2}(2,1,0,\ldots,1,1,1), \] \[R_{0,n-1}(2,0,1,\ldots,1,1,1),
Q_{0,1}(1,2,1,1,\ldots,1,1,0),Q_{0,2}(1,1,2,1,\ldots,1,1,0), \] \[ \ldots \ldots,\
Q_{0,n-3}(1,1,1,1,\ldots,2,1,0), Q_{0,n-2}(1,1,1,1,\ldots,1,2,0)\}\subset{{\bf N}^{n}}.\]
We will denote  by $A_{1}\in\ M_{n-1, n}({\bf R})$ the matrix with rows the coordinates of
points $ \{ \ R_{0,1}, R_{0,2},\ldots, R_{0,n-1} \}$ and for $2\leq \ i \leq \ n-1 $,
$A_{i}\in\ M_{n-1, n}({\bf R})$ the matrix with rows the coordinates of points $
\{ \ R_{0,1}, \ldots, R_{0,n-i},Q_{0,1}, Q_{0,2}, \ldots, Q_{0,i-1} \}$, that is:
\[ A_{1} = \left(
         \begin{array}{ccccccccc}
           2 & 1 & 1 & 1 & \ldots & \ldots & 1 & 1 & 0 \\
           2 & 1 & 1 & 1 & \ldots & \ldots & 1 & 0 & 1 \\
           2 & 1 & 1 & 1 & \ldots & \ldots & 0 & 1 & 1 \\
           . & . & . & . & . & . & . & . & . \\
           2 & 1 & 1 & 0 & \ldots  & \ldots & 1 & 1 & 1 \\
           2 & 1 & 0 & 1 & \ldots  & \ldots & 1 & 1 & 1 \\
           2 & 0 & 1 & 1 & \ldots  & \ldots & 1 & 1 & 1 \\
         \end{array}
       \right) \]
       and for $2\leq \ i \leq \ n-1 $\\

\[
\begin{array}{cccccc}
      &  &  &  & \downarrow(n-i)^{th} column &
\end{array}
\]
\[ A_{i} = \left(
   \begin{array}{ccccccccccc}
     2 & 1 & 1 &  \ldots & 1 & 1 & 1 \ \ldots  \ \ 1 & 1 & 0  \\
     2 & 1 & 1 &  \ldots & 1 & 1 & 1 \ \ldots  \ \ 1 & 0 & 1  \\
     2 & 1 & 1 &  \ldots & 1 & 1 & 1 \ \ldots  \ \ 0 & 1 & 1  \\
     . & . & . &  \ldots & . & . & . \ \ldots  \ \ . & . & .  \\
     2 & 1 & 1 &  \ldots & 1 & 1 & 0 \ \ldots  \ \ 1 & 1 & 1 \\
     1 & 2 & 1 &  \ldots & 1 & 1 & 1 \ \ldots  \ \ 1 & 1 & 0 \\
     1 & 1 & 2 &  \ldots & 1 & 1 & 1 \ \ldots  \ \ 1 & 1 & 0 \\
     . & . & . &  \ldots & . & . & . \ \ldots  \ \ . & . & .  \\
     1 & 1 & 1 &  \ldots & 2 & 1 & 1 \ \ldots  \ \ 1 & 1 & 0 \\
     1 & 1 & 1 &  \ldots & 1 & 2 & 1 \ \ldots  \ \ 1 & 1 & 0 \\
   \end{array}
 \right)\begin{array}{c}
           \\
           \\
 \leftarrow\ (n-i+1)^{th}\ line \\
           \\
        \end{array}
\]
Let $T_{i}$ be the linear transformation from ${\bf R^{n}}$ into ${\bf R^{n-1}}$ defined
by $T_{i}(x)=A_{i}x$ for all $1\leq \ i \leq \ n-1 $.

Let $i,j \in \bf N$, $1\leq i,j\leq n$ we denote by $e_{i,j}$ the
matrix from $M_{n-1}({\bf R})$ with the entries: one, on position
$"i, j"$  and zero in else; we denote by $T_{i,j}(a)$ the matrix
\[ T_{i,j}(a)=I_{n-1}+ae_{i,j} \in M_{n-1}(\bf R).\]
By $P_{i,j}$ we denote the matrix from $ M_{n-1}(\bf R)$ defined by
\[P_{i,j}=I_{n-1}-e_{i,i}-e_{j,j}+e_{i,j}+e_{j,i}\]
\begin{lemma}
$a)$ The set of points $ \{ \ R_{0,1}, \ldots, R_{0,n-i},Q_{0,1},
Q_{0,2}, \ldots, Q_{0,i-1} \}$ for all $2\leq \ i \leq \ n-1 $ and
$ \{ \ R_{0,1}, R_{0,2},\ldots, R_{0,n-1} \}$ are linearly
independent.\\ $b)$ For $1\leq \ i \leq \ n-1 $ the equations of
the hyperplanes generated by the points \\ $ \{ \ O, R_{0,1},
R_{0,2} \ldots, R_{0,n-i},Q_{0,1}, Q_{0,2}, \ldots, Q_{0,i-1} \}$
are :
\[ H_{[i]}:= -(n-i-1)\sum_{j=1}^{i}x_{j} + (i+1)\sum_{j=i+1}^{n}x_{j}=0,\]
where $[i]\ is\ the\ set\ [i]:=\{1,\ldots \,i\}. $

\end{lemma}
\begin{proof}
$a)$ The set of points are linearly independent iff the matrices
with rows the coordinates of the points have the rank $n-1 .$\\
Using elementary row operations to the matrix $A_{1}$ we have:\\ $
B_{1}=U_{1}A_{1},$ where $U_{1}\in\ M_{n-1}({\bf R})$ \[ U_{1}=
\prod_{2\leq i\leq \lfloor\frac{n}{2}\rfloor}
P_{i,n-i+1}\prod_{i=2}^{n-1}T_{n-i+1,1}(-1),\] and $ \lfloor c
\rfloor$ is the greatest integer $\leq c.$  So $B_{1}$ is :

\[ B_{1} = \left(
         \begin{array}{ccccccccc}
           2 & 1 & 1 & 1 & \ldots & \ldots & 1 & 1 & 0 \\
           0 & -1 & 0 & 0 & \ldots & \ldots & 0 & 0 & 1 \\
           0 & 0 & -1 & 0 & \ldots & \ldots & 0 & 0 & 1 \\
           . & . & . & . & . & . & . & . & . \\
           0 & 0 & 0 & 0 & \ldots  & \ldots & 0 & 0 & 1 \\
           0 & 0 & 0 & 0 & \ldots  & \ldots & -1 & 0 & 1 \\
           0 & 0 & 0 & 0 & \ldots  & \ldots & 0 & -1 & 1 \\
         \end{array}
       \right) \]

\[\]

For $2\leq \ i \leq \ n-1,$ using elementary row operations to the matrix $A_{i}$ we have:
$ B_{i}=U_{i}A_{i},$ where $U_{i}\in\ M_{n-1}({\bf R})$,
\[ U_{i}=[ \prod_{j=i}^{n-2}( \prod_{k=1}^{i-1})P_{n-j+k-1,n-j+k}][ \prod_{k=2}^{i-1}
(
\prod_{j=n-i+k}^{n-1}T_{j,n-i+k-1}(-\frac{1}{k+1}))]\cdot\]\[\cdot(
\prod_{j=n-i+1}^{n-1} T_{j,1}(-\frac{1}{2}))(
\prod_{j=1}^{n-i}T_{j,1}(-1)), \] and so $B_{i}$ is :

\[
\begin{array}{cccccccccccc}
      &  &  &  &  &  &  &  &  &  &\downarrow(i+1)^{th} column &
\end{array}
\]

\[
 B_{i} = \left(
   \begin{array}{cccccccccccc}
     2 & 1 & 1 & 1 & \ldots & 1 & 1 & 1 \ \ldots  \ \ 1 & 1 & 0  \\
     0 & \frac{3}{2} & \frac{1}{2} & \frac{1}{2} & \ldots & \frac{1}{2} & \frac{1}{2} & \frac{1}{2} \ \ldots  \ \ \frac{1}{2} & \frac{1}{2} & 0  \\
     0 & 0 & \frac{4}{3} & \frac{1}{3} & \ldots & \frac{1}{3} & \frac{1}{3} & \frac{1}{3} \ \ldots  \ \ \frac{1}{3} & \frac{1}{3} & 0  \\
     . & . & . & . & \ldots & . & . & \ \ldots  \ \ & . & .  \\
     0 & 0 & 0 & 0 & \ldots & \frac{i}{i-1} & \frac{1}{i-1} & \frac{1}{i-1} \ldots  \ \frac{1}{i-1} & \frac{1}{i-1} & 0 \\
     0 & 0 & 0 & 0 & \ldots & 0 & \frac{i+1}{i} &\frac{1}{i} \ \ldots \ \frac{1}{i} & \frac{1}{i} & 0 \\
     0 & 0 & 0 & 0 & \ldots & 0 & 0 & -1 \ \ldots \ 0 & 0 & 1 \\
     . & . & . & . & \ldots & . & . & . \ \ldots  \ \ . & . & .  \\
     0 & 0 & 0 & 0 & \ldots & 0 & 0 & \ \ \ \ 0 \ \ \ldots  -1 & 0 & 1 \\
     0 & 0 & 0 & 0 & \ldots & 0 & 0 & \ 0 \ \ldots  \ 0 & -1 & 1 \\
   \end{array}
 \right)\begin{array}{c}
           \\
           \\
 \leftarrow\ (i)^{th}\ line \\
           \\
        \end{array}
\]

\[\]

Since the rank of $B_{i}$ is $n-1$ then the rank of $A_{i}$ is
$n-1$ for all $1\leq \ i \leq \ n-1 .$

$b)$ The hyperplane generated by the points $ \{ \ R_{0,1},
\ldots, R_{0,n-i},Q_{0,1}, Q_{0,2}, \ldots, Q_{0,i-1} \}$ has the
normal vector the generator of the subspace $Ker(T_{i})$.

For $1\leq \ i \leq \ n-1 $ using $a)$ we obtain that
\[ Ker(T_{i})=\{ x\in{\bf R^{n}}|T_{i}(x)=0 \}=\{ x\in{\bf R^{n}}|A_{i}x=0 \}=\{ x\in{\bf R^{n}}
|B_{i}x=0 \}\] that is \[ x_{n}=x_{n-1}=\ldots =x_{i+1}=(i+1)\alpha \] and \[ x_{i}=x_{i-1}=
\ldots =x_{1}=-(n-i-1)\alpha ,\] where $\alpha \in \bf R.$

Thus for $1\leq \ i \leq \ n-1 $ the equations of hyperplanes generated by the points
$ \{ \ R_{0,1}, \ldots, R_{0,n-i},Q_{0,1}, Q_{0,2}, \ldots, Q_{0,i-1} \}$ are :
\[ H_{[i]}:= -(n-i-1)\sum_{j=1}^{i}x_{j} + (i+1)\sum_{j=i+1}^{n}x_{j}=0,\]

\end{proof}
For $1\leq k \leq n-1$ we define two types  of  sets of points:

$1)$ \[\{ \ R_{k,1}, R_{k,2},\ldots, R_{k,n-1} \}\] is the set of
points whose the coordinates are the rows of the matrix
$A_{1}P_{1k+1}.$

$2)$ \[\{ \ Q_{k,1}, Q_{k,2},\ldots, Q_{k,n-2} \}\] is the set of
points whose the coordinates are the rows of the matrix $QM^{k}$,
where $M$ is the matrix  \[M\in\ M_{n}({\bf R}), M=
\prod_{i=1}^{n-1}P_{n-i,n-i+1} \] and $Q \in\ M_{n-2, n}({\bf R})$
is the matrix with rows the coordinates of points $ \{ \ Q_{1},
Q_{2}, \ldots, Q_{n-2} \}$.

For every $1\leq i \leq n-1$ we will denote by $\nu_{[i]}$ the normal of the hyperplane
$H_{[i]}$:
\[\begin{array}{cccccc}
     &  &  &  & \downarrow i^{th} column &
  \end{array}
\]
\[\nu_{[i]} = \left( \begin{array}{cccccc}
     -(n-i-1), & \ldots & ,-(n-i-1), & (i+1), & \ldots & ,(i+1)  \\
   \end{array}
   \right) \in {\bf R^{n}}
\]

For $i=1$ we will denote by $H_{\{k+1 \}}$ the hyperplane which
has the normal :
\[ \nu_{\{k+1 \}}:=\nu_{[i]}P_{1,k+1}=\nu_{[1]}P_{1,k+1}\] for all $1\leq k \leq n-1.$

For $2\leq i \leq n-1$ and $1\leq k \leq n-1$ we will denote by
$H_{\{\sigma^{k}(1), \sigma^{k}(2),\ldots ,\sigma^{k}(i) \}}$ the
hyperplane which has the normal :
\[\nu_{\{\sigma^{k}(1), \sigma^{k}(2),\ldots ,\sigma^{k}(i) \}}:=\nu_{[i]}M^{k},\]
where $\sigma \in S_{n}$ is the product of transposition :\[\sigma :=\prod_{i=1}^{n-1}(i,i+1).\]

\begin{lemma}
$a)$ For $1\leq k \leq n-1$ and $2\leq \ i \leq \ n-1 $ the set of
points $ \{ \ R_{k,1}, \ldots, R_{k,n-i},\\Q_{k,1}, Q_{k,2},
\ldots, Q_{k,i-1} \}$  and $ \{ \ R_{k,1}, R_{k,2},\ldots,
R_{k,n-1} \}$ are linearly independent.\\ $b)$ For $1\leq k \leq
n-1$ and $2\leq \ i \leq \ n-1 $ the equation of hyperplane
generated by the points $ \{ \ O, R_{k,1}, R_{k,2} \ldots,
R_{k,n-i},Q_{k,1}, Q_{k,2}, \ldots, Q_{k,i-1} \}$ is :
\[ H_{\{\sigma^{k}(1), \sigma^{k}(2),\ldots ,\sigma^{k}(i) \}} := <\nu_{\{\sigma^{k}(1),
\sigma^{k}(2),\ldots ,\sigma^{k}(i) \}},x>=0\] where $O$ is zero
point, $O(0, 0, \ldots , 0)$ and $\sigma \in S_{n}$ is the product
of transposition :\[\sigma :=\prod_{i=1}^{n-1}(i,i+1).\] For
$1\leq k \leq n-1$ the equation of hyperplane generated by the
points \\ $ \{ \ O, R_{k,1}, R_{k,2} \ldots, R_{k,n-1}\}$ is
:\[H_{\{k+1 \}}:=<\nu_{\{k+1 \}},x>=0\]

\end{lemma}

\begin{proof}
$a)$ Since for any $1\leq k \leq n-1$ the matrix $P_{1,k+1}$ ,
$M^{k}$ are invertible and the set of points $ \{ \ R_{0,1},
\ldots, R_{0,n-i},Q_{0,1}, Q_{0,2}, \ldots, Q_{0,i-1} \}$ , $ \{ \
R_{0,1}, R_{0,2},\ldots, R_{0,n-1} \}$ are linearly independent
then the set of points $ \{ \ R_{k,1}, \ldots, R_{k,n-i},Q_{k,1},
Q_{k,2}, \ldots, Q_{k,i-1} \}$  and \\ $ \{ \ R_{k,1},
R_{k,2},\ldots, R_{k,n-1} \}$ are linearly independent.

$b)$ Since for any $1\leq k \leq n-1$ and $2\leq \ i \leq \ n-1 $
the matrix $M^{k}$ are invertible then the hyperplane generated by
the points $ \{ \ O, R_{k,1}, \ldots, R_{k,n-i},Q_{k,1}, \ldots,
Q_{k,i-1} \}$ has the normal vector obtained by normal vector
$\nu_{[k]}$ multiple to the right with $M^{k}$. For any $1\leq k
\leq n-1$ the matrix $P_{1,k+1}$ are invertible, then the
hyperplane generated by the points $ \{ \ O, R_{k,1}, R_{k,2}
\ldots, R_{k,n-1}\}$ has the normal vector obtained by normal
vector $\nu_{[1]}$ multiple to the right with $P_{1,k+1}$.
\end{proof}

\begin{lemma}
Any point $P \in {\bf N^{n}}$, $n\geq 3$ which lies on the
hyperplane $H : x_{1}+x_{2}+ \ldots + x_{n}-n=0$ such that its
coordinates are in the set $\{0,1,2\}$ and has at least one
coordinate equal to two lies on the hyperplane $H_{\{k\}}=0$ for a
integer $k \in \{1,2,...,n\}.$
\end{lemma}

\begin{proof}
Let $k \in \{1,2,...,n\}$ be the first position of $"2"$ that
appears in the coordinate of a point $P \in {\bf N^{n}}$. Since
the equation of the hyperplane $H_{\{k\}}$ is: \[ H_{\{k\}}=
\sum_{i=1}^{k-1}2x_{i}-(n-2)x_{k}+ \sum_{i=k+1}^{n}2x_{i}=0 \] it
results that
\[-2(n-2)+2\sum_{i=1,i\neq k}{n}a_{i}=-2(n-2)+2(n-2)=0,\] where $P=(a_{1},a_{2},\ldots,a_{n})
\in H \ with \ a_{i}\in \{0,1,2\}$ and which has at least one
coordinate equal to two.
\end{proof}
\section{The main result}
First let us fix some notations that will be used throughout the
remaining of this note. Let $K$ be field and
$K[x_{1},x_{2},\ldots,x_{n}]$ a polynomial ring with coefficients
in $K$. Let $n\geq 2$ a positive integer and ${\bf{\mathcal{A}}}$
the collection of sets:
\[{{\bf{\mathcal{A}}}}=\{\{1,2\},\{2,3\},\ldots,\{n-1,n\},\{n,1\}\}. \]
We will denote by $K[{\bf{\mathcal{A}}}]$ the $K-$algebra generated by
$x_{i_{1}}x_{i_{2}}...x_{i_{n}}$
with $i_{1}\in\{1,2\},i_{2}\in\{2,3\},\ldots,i_{n-1}\in\{n-1,n\}, i_{n}\in\{1,n\}. $
This $K-$algebra represent the base ring associated to transversal polymatroid presented by
${\bf{\mathcal{A}}}.$

Given $A\in {\bf{N^{n}}}$ finite, we define $C_{A}$ to be the subsemigroup of ${\bf{N^{n}}}$
generated by $A:$
\[C_{A}=\sum_{\alpha \in A}{\bf{N}}\alpha\] thus the $cone$ generated by $C_{A}$ is:
\[{\bf{R_{+}}}C_{A}={\bf{R_{+}}}A=\{ \sum a_{i}\gamma_{i}\ | \ a_{i}\in {\bf R_{+}},
\gamma_{i}\in A\}.\]

With this notation we state our main result:

\begin{theorem}
Let $A=\{ log(x_{i_{1}}x_{i_{2}}...x_{i_{n}}) \ | \ i_{1}\in\{1,2\},i_{2}\in\{2,3\},\ldots,
i_{n-1}\in\{n-1,n\}, i_{n}\in\{1,n\}\}\subset {\bf{N^{n}}}$ the exponent set of the generators
of $K-$algebra $K[{\bf{\mathcal{A}}}]$ and $N=\{\nu_{\{k+1 \}}, \nu_{\{\sigma^{k}(1),
\sigma^{k}(2),\ldots ,\sigma^{k}(i) \}} \ | \ 0\leq k \leq n-1, 2\leq i \leq n-1\}$, then
\[{\bf{R_{+}}}C_{A}= \bigcap_{a\in N}H^{+}_{a},\]
such that $H^{+}_{a}$ with $a\in N$ are just the facets of the
cone ${\bf{R_{+}}}C_{A}.$
\end{theorem}
\begin{proof}
Since $A=\{ log(x_{i_{1}}x_{i_{2}}...x_{i_{n}}) \ | \ i_{1}\in\{1,2\}, \ i_{2}\in\{2,3\}\ ,
\ldots, \ i_{n-1}\in\{n-1,n\}, \  i_{n}\in\{1,n\}\}\subset {\bf{N^{n}}}$ is the exponent set
of the generators of $K-$algebra $K[{\bf{\mathcal{A}}}]$, then the set $ \{ \ R_{0,1}, R_{0,2},
\ldots, R_{0,n-2}, R_{0,n-1}, I\}\subset A$, where $I(1,1,\ldots,1)\in {\bf{N^n}}.$

{\bf{First step.}}

We must show that the dimension of the cone ${\bf{R_{+}}}C_{A}$ is $dim({\bf{R_{+}}}C_{A})=n.$\\
We will denote by $\widetilde{A}\in  M_{n}({\bf R})$ the matrix with rows the coordinates of
points $ \{ \ R_{0,1}, R_{0,2},$ $\ldots, R_{0,n-2}, R_{0,n-1}, I\}.$
Using elementary row operations to the matrix $\widetilde{A}$ we have:
$\widetilde{B}=\widetilde{U}\widetilde{A}$, where $\widetilde{U}\in\ M_{n}({\bf R})$ is
invertible matrix: \[ \widetilde{U}= (\prod_{i=2}^{n-1}T_{n-i+1,1}(-1)))(T_{n,1}(-\frac{1}{2}))
(\prod_{2\leq i\leq \lfloor\frac{n}{2}\rfloor} P_{i,n-i+1})(\prod_{i=2}^{n-1}T_{n,n-i+1}
(\frac{1}{2})),\] where $ \lfloor c \rfloor$ is the greatest integer $\leq c.$

So $\widetilde{B}$ is:

\[ \widetilde{B} = \left(
         \begin{array}{ccccccccc}
           2 & 1 & 1 & 1 & \ldots & \ldots & 1 & 1 & 0 \\
           0 & -1 & 0 & 0 & \ldots & \ldots & 0 & 0 & 1 \\
           0 & 0 & -1 & 0 & \ldots & \ldots & 0 & 0 & 1 \\
           . & . & . & . & . & . & . & . & . \\
           0 & 0 & 0 & 0 & \ldots  & \ldots & 0 & 0 & 1 \\
           0 & 0 & 0 & 0 & \ldots  & \ldots & -1 & 0 & 1 \\
           0 & 0 & 0 & 0 & \ldots  & \ldots & 0 & -1 & 1 \\
           0 & 0 & 0 & 0 & \ldots  & \ldots & 0 & 0 & \frac{n}{2} \\
         \end{array}
       \right) \]
Then the dimension of the cone ${\bf{R_{+}}}C_{A}$ is:
\[dim({\bf{R_{+}}}C_{A})=
rank(\widetilde{A})=rank(\widetilde{B})=n\] since
$det(\widetilde{B})=(-1)^{n}n.$\\

{\bf{Second step.}}

We must show that $H_{a}\cap{\bf{R_{+}}}C_{A}$ with $a\in N$ are
precisely the facets of the cone ${\bf{R_{+}}}C_{A}.$ This is
equivalent with to show that ${\bf{R_{+}}}C_{A}\subset H^{+}_{a}$
and $dim H_{a}\cap {\bf{R_{+}}}C_{A}=n-1$ $\forall$ $a\in N.$\\The
fact that $dim H_{a}\cap {\bf{R_{+}}}C_{A}=n-1$ $\forall$ $a\in N$
it is clear from Lemma 3.1 and Lemma 3.2.

For $1\leq k\leq \lfloor\frac{n}{2}\rfloor$ and $1\leq i_{1}<i_{2}< \ldots < i_{2k-1}<
i_{2k}\leq n$
let \[I_{i_{1}i_{2}\ldots i_{2k-1}i_{2k}}=I + (e_{i_{1}}-e_{i_{2}}) + (e_{i_{3}}-e_{i_{4}})+
\ldots + (e_{i_{2k-1}}-e_{i_{2k}})\] and \[I^{'}_{i_{1}i_{2}\ldots i_{2k-1}i_{2k}}=
I + (e_{i_{2}}-e_{i_{1}}) + (e_{i_{4}}-e_{i_{3}})+\ldots + (e_{i_{2k}}-e_{i_{2k-1}}),\]
where $I=I(1,1,\ldots,1)\in {\bf{N^n}}$ and $e_{i}$ is the $ith$ unit vector.

We set $A^{'}=$\[\{I, I_{i_{1}i_{2}\ldots i_{2k-1}i_{2k}}, I^{'}_{i_{1}i_{2}
\ldots i_{2k-1}i_{2k}} | 1\leq k\leq \lfloor\frac{n}{2}\rfloor \ and\ 1\leq i_{1}<i_{2}<
\ldots < i_{2k-1}< i_{2k}\leq n \}.\]

We claim that $A=A^{'}.$

Let \[m_{i_{1}i_{2}\ldots
i_{2k-1}i_{2k}}=\prod_{s=1}^{\lfloor\frac{n}{2}\rfloor}m_{s}, \
m^{'}_{i_{1}i_{2}\ldots
i_{2k-1}i_{2k}}=\prod_{s=1}^{\lfloor\frac{n}{2}\rfloor}m_{s}^{'}\]
where \[ m_{s}=x_{k_{i_{2s-2}+1}}\ldots
x_{k_{i_{2s-1}-2}}x_{i_{2s-1}}^2x_{k_{i_{2s-1}+1}}\ldots
x_{k_{i_{2s}-2}}x_{i_{2s}-1}x_{i_{2s}+1},\] \[
m_{s}^{'}=x_{k_{i_{2s-2}+1}}\ldots
x_{k_{i_{2s-1}-2}}x_{i_{2s-1}-1}x_{i_{2s-1}+1}x_{k_{i_{2s-1}+1}}\ldots
x_{k_{i_{2s}-2}}x_{i_{2s}}^{2}\] for all $1\leq k, s\leq
\lfloor\frac{n}{2}\rfloor,$ $i_{0}=0$ and $k_{j}\in \{j,j+1\}$ for
$1\leq j\leq n.$ Evidently $log(m_{i_{1}i_{2}\ldots
i_{2k-1}i_{2k}}),\\ log(m^{'}_{i_{1}i_{2}\ldots
i_{2k-1}i_{2k}})\in A.$

Since $log(m_{i_{1}i_{2}\ldots i_{2k-1}i_{2k}})=I_{i_{1}i_{2}\ldots i_{2k-1}i_{2k}}$ and
$log(m^{'}_{i_{1}i_{2}\ldots i_{2k-1}i_{2k}})=I^{'}_{i_{1}i_{2}\ldots i_{2k-1}i_{2k}}$
for all $1\leq k\leq \lfloor\frac{n}{2}\rfloor$ \ and\ $1\leq i_{1}<i_{2}< \ldots < i_{2k-1}<
i_{2k}\leq n,$  then $A^{'}\subset A.$

But the cardinal of $A$ is $\sharp(A)=2^{n}-1$ and since
\[\sum_{s=1}^{\lfloor\frac{n}{2}\rfloor}\binom{n}{2s}=2^{n-1}-1,\] the cardinal of $A^{'}$is:
\[\sharp(A^{'})=1+2\sum_{s=1}^{\lfloor\frac{n}{2}\rfloor}\binom{n}{2s}=2^{n}-1.\]
Thus $A=A^{'}.$

Now we start to prove that ${\bf{R_{+}}}C_{A}\subset H^{+}_{a}$
for all $a\in N.$\\Observe that \[<\nu_{\{p+1\}},I>\ =\
<\nu_{\{\sigma^{p}(1), \sigma^{p}(2),\ldots ,\sigma^{p}(i) \}},I>\
= n>0,\]
for any $0\leq p \leq n-1, $ $1\leq i \leq n-1. $\\
Let $0\leq p \leq n-1.$ We claim that:\[<\nu_{\{p+1\}},I_{i_{1}i_{2}\ldots i_{2k-1}i_{2k}}> \
\geq 0 \ and \ <\nu_{\{p+1\}},I^{'}_{i_{1}i_{2}\ldots i_{2k-1}i_{2k}}> \ \geq 0, \] for any
$1\leq k\leq \lfloor\frac{n}{2}\rfloor$ and
$1\leq i_{1}<i_{2}< \ldots < i_{2k-1}< i_{2k}\leq n.$ We prove the first inequality.
The proof of the second inequality is the same. \\We have three possibilities:\\ $1)$
If $<I_{i_{1}i_{2}\ldots i_{2k-1}i_{2k}}, e_{p+1}>=0,$
then $<\nu_{\{p+1\}},I_{i_{1}i_{2}\ldots i_{2k-1}i_{2k}}>=2n>0.$\\$2)$
If $<I_{i_{1}i_{2}\ldots i_{2k-1}i_{2k}}, e_{p+1}>=1,$
then $<\nu_{\{p+1\}},I_{i_{1}i_{2}\ldots i_{2k-1}i_{2k}}>=n>0.$\\$3)$
If $<I_{i_{1}i_{2}\ldots i_{2k-1}i_{2k}}, e_{p+1}>=2,$
then $<\nu_{\{p+1\}},I_{i_{1}i_{2}\ldots i_{2k-1}i_{2k}}>=0.$

Let $0\leq p \leq n-1$ and $2\leq i \leq n-1$ fixed.

We claim that:
\[<\nu_{\{\sigma^{p}(1), \sigma^{p}(2),\ldots ,\sigma^{p}(i) \}}, I_{i_{1}i_{2}\ldots
i_{2k-1}i_{2k}}> \\ \geq 0 \]and \[<\nu_{\{\sigma^{p}(1),
\sigma^{p}(2),\ldots ,\sigma^{p}(i) \}}, I^{'}_{i_{1}i_{2}\ldots i_{2k-1}i_{2k}}> \\ \geq 0 \]
for any $1\leq k\leq \lfloor\frac{n}{2}\rfloor$ and
$1\leq i_{1}<i_{2}< \ldots < i_{2k-1}< i_{2k}\leq n.$\\ We prove the first inequality.
The proof of the second inequalities is the same.

We have:
\[ <\nu_{\{\sigma^{p}(1), \sigma^{p}(2),\ldots ,\sigma^{p}(i) \}},
I_{i_{1}i_{2}\ldots i_{2k-1}i_{2k}}>=H_{\{\sigma^{p}(1), \sigma^{p}(2),\ldots ,\sigma^{p}(i) \}}
(I_{i_{1}i_{2}\ldots i_{2k-1}i_{2k}})=\]
\[-(n-i-1)\sum_{s=1}^{i}<I_{i_{1}i_{2}\ldots i_{2k-1}i_{2k}}, e_{\sigma^{p}(s)}>+ \ (i+1)
\sum_{s=i+1}^{n}<I_{i_{1}i_{2}\ldots i_{2k-1}i_{2k}}, e_{\sigma^{p}(s)}>.\]

Let \[ \Gamma=\{s | <I_{i_{1}i_{2}\ldots i_{2k-1}i_{2k}},
e_{\sigma^{p}(s)}>=2, 1\leq s \leq i\}\] the set of indices of
$I_{i_{1}i_{2}\ldots i_{2k-1}i_{2k}}$ where the coordinates are
equal to two.\\
If the cardinal of $\Gamma$ is zero, then there exists at most
one index $i_{2t-1}\in \{\sigma^{p}(1), \sigma^{p}(2),\ldots\\
,\sigma^{p}(i)\}$ with $1\leq t\leq \lfloor\frac{n}{2}\rfloor.$
Else we have two possibilities:
\\$1)$There exists at least two indices
$i_{2t-1}, i_{2t_{1}-1}\in \{\sigma^{p}(1), \sigma^{p}(2),\ldots
,\sigma^{p}(i)\}$ with $1\leq t < t_{1}\leq
\lfloor\frac{n}{2}\rfloor$ and since $\sigma^{p}(s)=(p+s) \ mod \
n,$ then there exists $1\leq t_{2}\leq \lfloor\frac{n}{2}\rfloor$
such that $i_{2t_{2}}\in \{\sigma^{p}(1), \sigma^{p}(2),\ldots
,\sigma^{p}(i)\}$ and thus $<I_{i_{1}i_{2}\ldots i_{2k-1}i_{2k}},
e_{\sigma^{p}(i_{2t_{2}})}> =2,$ which it is false.
\\$2)$There exists at least two indices $i_{2t-1}, i_{2t_{1}}\in \{\sigma^{p}(1),
\sigma^{p}(2),\ldots ,\sigma^{p}(i)\}$ with $1\leq t, \ t_{1}\leq
\lfloor\frac{n}{2}\rfloor.$ Then like in the first case we have
$<I_{i_{1}i_{2}\ldots i_{2k-1}i_{2k}}, e_{\sigma^{p}(i_{2t_{1}})}>
=2,$ which it is false.
\\When for any $1\leq k\leq \lfloor\frac{n}{2}\rfloor,$ $i_{2k-1} \not \in \{\sigma^{p}(1),
\sigma^{p}(2),\ldots ,\sigma^{p}(i)\},$ then \[ <\nu_{\{\sigma^{p}(1), \sigma^{p}(2),\ldots ,
\sigma^{p}(i) \}}, I_{i_{1}i_{2}\ldots i_{2k-1}i_{2k}}>=-(n-i-1)i+(i+1)(n-i)=n>0.\]\\
When there exists just one index $i_{2t-1}\in \{\sigma^{p}(1),
\sigma^{p}(2),\ldots , \sigma^{p}(i)\}$ with $1\leq t\leq
\lfloor\frac{n}{2}\rfloor,$ then \[ <\nu_{\{\sigma^{p}(1),
\sigma^{p}(2),\ldots ,\sigma^{p}(i) \}}, I_{i_{1}i_{2}\ldots
i_{2k-1}i_{2k}}>=-(n-i-1)(i-1)+(i+1)(n-i+1)=2n>0.\]

 If the cardinal of $\Gamma,$ $\sharp(\Gamma)=t\geq 1,$ then we have two possibilities:\\
$1)$ If $\{1\leq i_{1}<i_{2}<\ldots
<i_{2t-3}<i_{2t-2}<i_{2t-1}\leq n\}\subset \{\sigma^{p}(1),
\sigma^{p}(2),\ldots ,\sigma^{p}(i)\},$ \\ then we have: \[
<\nu_{\{\sigma^{p}(1), \sigma^{p}(2),
\ldots ,\sigma^{p}(i) \}}, I_{i_{1}i_{2}\ldots i_{2k-1}i_{2k}}>=-(n-i-1)(i+1)+(i+1)(n-i-1)=0.\]\\
$2)$ If $\{1\leq i_{1}<i_{2}<\ldots <i_{2t-1}<i_{2t}\leq n\}\subset \{\sigma^{p}(1),
\sigma^{p}(2),\ldots ,\sigma^{p}(i)\},$ then we have: \[ <\nu_{\{\sigma^{p}(1),
\sigma^{p}(2),\ldots ,\sigma^{p}(i) \}}, I_{i_{1}i_{2}\ldots i_{2k-1}i_{2k}}>=
-(n-i-1)(i)+(i+1)(n-i)=n>0.\]

Thus we have: \[{\bf{R_{+}}}C_{A}\subset \bigcap_{a\in
N}H^{+}_{a}\]

Finally let us prove the converse inclusion.\\ This is equivalent
with the fact that the extremal rays of the cone \[\bigcap_{a\in
N}H^{+}_{a}\] are in ${\bf{R_{+}}}C_{A}.$

Let $1\leq k\leq \lfloor\frac{n}{2}\rfloor,$  $1\leq i_{1}<i_{2}<
\ldots < i_{2k-1}< i_{2k}\leq n$ and we consider the \\following
hyperplanes:\\ $a)$ $H_{\{[i_{2s-1}]\backslash [j]\}}$ if
$j\in\{i_{2s-2},\ldots \,i_{2s-1}-1\}$ and $1\leq s \leq k ,$\\
$b)$ $H_{\{[j]\backslash [i_{2s-1}-1]\}}$ if
$j\in\{i_{2s-1}+1,\ldots ,i_{2s}-1\}$ and $1\leq s \leq k ,$\\
$c)$ $H_{\{[i_{2k-1}]\cup ( [n] \backslash [j])\}}$ if
$j\in\{i_{2k},\ldots \,n-1\} ,$\\ $d)$ $H_{\{i_{2s}\}}$ for $1\leq
s \leq k-1 ;$ where $[i]:=\{1,\ldots , i\},$ $i_{0}=0$ and
$[0]=\emptyset.$\\
We claim that the point $I_{i_{1}i_{2}\ldots i_{2k-1}i_{2k}}$
belong to this hyperplanes.\\
$a)$ Let $j\in\{i_{2s-2},\ldots \,i_{2s-1}-1\}$ and $1\leq s \leq
k ,$ then \[<H_{\{[i_{2s-1}]\backslash [j]\}} ,
I_{i_{1}i_{2}\ldots i_{2k-1}i_{2k}}>=\ <H_{\{j+1,\ldots
,i_{2s-1}\}} , I_{i_{1}i_{2}\ldots
i_{2k-1}i_{2k}}>=<-(n-(i_{2s-1}-j)-1)\]\[\sum_{t\in\{j+1,\ldots ,
i_{2s-1} \}} x_{t}  \\ + (i_{2s-1}-j+1)\sum_{t\in [n]\backslash
\{j+1,\ldots , i_{2s-1} \}} x_{t}\ , \ I_{i_{1}i_{2}\ldots
i_{2k-1}i_{2k}}>=-(n-i_{2s-1}+j-1)(i_{2s-1}\] $-j+1)\ + \
(i_{2s-1}-j+1)(n-(i_{2s-1}-j)+1)=0,$ since
\[
\begin{array}{ccccccccccccccccccccccccccc}
    &   &   &   &   &   &   &   &   &  &  &  &  &  &  \downarrow j+1^{th} &   &   &   &   & \downarrow i_{2s-1}^{th} &   &   &
\end{array}
\]
\[
\begin{array}{ccccccccccccccc}
  I_{i_{1}i_{2}\ldots i_{2k-1}i_{2k}} & = & ( & \ldots & , & 1 & , & \ldots & , & 1 & , & 2 &  , & \ldots &
  ).
\end{array}
\]
$b)$ Let $j\in\{i_{2s-1}+1,\ldots \,i_{2s}-1\}$ and $1\leq s \leq
k ,$ then \[<H_{\{[j]\backslash [i_{2s-1}-1]\}} ,
I_{i_{1}i_{2}\ldots i_{2k-1}i_{2k}}>=\ <H_{\{i_{2s-1},\ldots ,j\}}
, I_{i_{1}i_{2}\ldots
i_{2k-1}i_{2k}}>=<-(n-(j-i_{2s-1}+1)-1)\]\[\sum_{t\in\{i_{2s-1},\ldots
,j\}} x_{t}  \\ + (j-i_{2s-1}+1+1)\sum_{t\in [n]\backslash
\{i_{2s-1},\ldots ,j \}} x_{t}\ , \ I_{i_{1}i_{2}\ldots
i_{2k-1}i_{2k}}>=-(n-(j-i_{2s-1}+1)-1)(j- \]$-i_{2s-1}+1+1)\ + \
(j-i_{2s-1}+1+1)(n-(j-i_{2s-1}+1+1))=0,$ since
\[
\begin{array}{ccccccccccccccccccccccccccc}
    &   &   &   &   &   &   &   &   &  &  &  &  &  \downarrow i_{2s-1}^{th} &   &   &   &   &   &\downarrow j^{th} &   &   &
\end{array}
\]
\[
\begin{array}{ccccccccccccccc}
  I_{i_{1}i_{2}\ldots i_{2k-1}i_{2k}} & = & ( & \ldots & , & 2 & , & 1 & , & \ldots & , & 1 &  , & \ldots &
  ).
\end{array}
\]
$c)$ Let $j\in\{i_{2k},\ldots \,n-1\},$ then $<H_{\{[i_{2k-1}]\cup
( [n] \backslash [j])\}} , I_{i_{1}i_{2}\ldots
i_{2k-1}i_{2k}}>=<H_{\{1,\ldots, i_{2k-1},j+1,\ldots, n\}},$
\[I_{i_{1}i_{2}\ldots
i_{2k-1}i_{2k}}>=<-(n-(i_{2k-1}+n-j)-1)\sum_{t\in[i_{2k-1}]\cup([n]
\backslash [j])}x_{t}+(i_{2k-1}+n-j+1)\sum_{t \in \{i_{2k-1}+1,
\ldots, j\}}x_{t} \ ,\]\[ I_{i_{1}i_{2}\ldots
i_{2k-1}i_{2k}}>=-(j-i_{2k-1}-1)(i_{2k-1}+n-j+1)+(i_{2k-1}+n-j+1)(j-(i_{2k-1}+1)+1-1)=0,
\] since
\[
\begin{array}{ccccccccccccccccccccccccccccccccccccccccccccccccccccccccccccccccccccccccc}
    &   &   &   &   &   &   &   &  &  &  &  &  &  & \downarrow i_{2k-1}^{th} &   &   &   &   &   &  &  &  & \downarrow i_{2k}^{th} &   &   &
    &   &   &  &  \downarrow j+1^{th} &   &   &   &   &
\end{array}
\]
\[
\begin{array}{ccccccccccccccccccccccccc}
  I_{i_{1}i_{2}\ldots
i_{2k-1}i_{2k}} & = & ( & \ldots & , & 2 & , & 1 & , & \ldots & ,&
1 & , & 0 & , & 1 &
  , & \ldots & , & 1 & , & \ldots & , & 1 & ).
\end{array}
\]
$d)$ It is clear from Lemma 3.3.\\
Since the number of hyperplanes is
$\sum_{s=1}^{k}(i_{2s-1}-1-i_{2s-2}+1)+
\sum_{s=1}^{k}(i_{2s}-1-(i_{2s-1}+1)+1) + (n-1-i_{2k}+1)
+k-1=\sum_{s=1}^{k}(i_{2s}-i_{2s-2}) -k + n - i_{2k} + k -1 =
n-1,$ then
\[\bigcap_{s=1}^{k}(\bigcap_{j=i_{2s-2}}^{i_{2s-1}-1}(H_{\{[i_{2s-1}]\backslash [j]})\cap\bigcap_{j=i_{2s-1}+1}^{i_{2s}-1}(H_{\{[j]\backslash [i_{2s-1}-1]\}}))
\cap \bigcap_{j=i_{2k}}^{n-1}(H_{\{[i_{2k-1}]\cup ( [n] \backslash
[j])\}})\cap
\bigcap_{s=1}^{k-1}(H_{\{i_{2s}\}})=\]$OI_{i_{1}i_{2}\ldots
i_{2k-1}i_{2k}}$ is a extreme ray  of the cone $\bigcap_{a\in
N}H^{+}_{a}$. But $OI_{i_{1}i_{2}\ldots i_{2k-1}i_{2k}}\in
{\bf{R_{+}}}C_{A}.$ Thus $\bigcap_{a\in
N}H^{+}_{a}={\bf{R_{+}}}C_{A}.$

\end{proof}
For use bellow we recall that $K-$algebra $K[{\bf{\mathcal{A}}}]$  is a normal domain
according to{\cite{HH}}.
\begin{definition}
Let $R$ be a polynomial ring over a field $K$ and $F$ a finite set
of monomials in $R.$ A decomposition
\[K[F]=\bigoplus_{i=0}^{\infty}K[F]_{i} \] of the $K-$ vector
space $K[F]$ is an $admissible\ grading$ if $k[F]$ is a positively
graded $K-$ algebra with respect to this decomposition and each
component $K[F]_{i}$ has a finite $K-$ basis consisting of
monomials.
\end{definition}
\begin{theorem}{\bf(Danilov, Stanley)}
Let $R=K[x_{1},\ldots,x_{n}]$ be a polynomial ring over a field $K$ and $F$ a finite set
of monomials in $R.$
If $K[F]$ is normal, then the canonical module $\omega_{K[F]}$ of $K[F],$ with respect to
an arbitrary admissible grading, can be expressed as an ideal of $K[F]$ generated by monomials
\[\omega_{K[F]}=(\{x^{a}| \ a\in {\bf{N}}A\cap ri({\bf{R_{+}}}A)\}),\] where $A=log(F)$ and
$ri({\bf{R_{+}}}A)$ denotes the relative interior of
${\bf{R_{+}}}A.$
\end{theorem}
\begin{corollary}
The canonical module $\omega_{K[{\bf{\mathcal{A}}}]}$ of $K[{\bf{\mathcal{A}}}]$
is $\omega_{K[{\bf{\mathcal{A}}}] }=(x_{1}x_{2}\ldots x_{n}) K[{\bf{\mathcal{A}}}].$
Thus $K-$ algebra $K[{\bf{\mathcal{A}}}]$  is Gorenstein ring.
\end{corollary}
\begin{proof}
Since \[<\nu_{\{p+1\}},I>\ =\ <\nu_{\{\sigma^{p}(1),
\sigma^{p}(2),\ldots ,\sigma^{p}(i) \}},I>\ = n>0,\] for any
$0\leq p \leq n-1, $ $1\leq i \leq n-1$ and since for any
$I_{i_{1}i_{2}\ldots i_{2k-1}i_{2k}}, I^{'}_{i_{1}i_{2} \ldots
i_{2k-1}i_{2k}}$ there exist two hyperplanes $H_{a}, H_{a^{'}}$
with $a, a^{'} \in N$ such that \[<I_{i_{1}i_{2}\ldots
i_{2k-1}i_{2k}}, H_{a}>=<I^{'}_{i_{1}i_{2}\ldots i_{2k-1}i_{2k}},
H_{a^{'}}>=0,\]then $I\in ri({\bf{R_{+}}}C_{A}) \ is\ the\ only\
point\ in\ relative\ interior\ of\ cone\ {\bf{R_{+}}}C_{A}.$ Thus
the canonical module is generated by one generator
$\omega_{K[{\bf{\mathcal{A}}}]}=(x_{1}x_{2}\ldots x_{n})
K[{\bf{\mathcal{A}}}].$ Thus $K-$ algebra $K[{\bf{\mathcal{A}}}]$
is Gorenstein ring.

\end{proof}

\begin{conjecture}
Let $n\in{\bf N}, A_{i}\subseteqq [n] \ with \ 1\leq i\leq n \ and
\ {\bf{\mathcal{\widetilde{A}}}}=\{A_{1},A_{2},\ldots,A_{n}\}.$ We
denote by $A=\{ log(x_{i_{1}}x_{i_{2}}...x_{i_{n}}) \ | \
i_{1}\in\{1,2\},i_{2}\in\{2,3\},\ldots, i_{n-1}\in\{n-1,n\},
i_{n}\in\{1,n\}\},$ $N=\{\nu_{\{k+1 \}}, \nu_{\{\sigma^{k}(1),
\sigma^{k}(2),\ldots ,\sigma^{k}(i) \}} \ | \ 0\leq k \leq n-1,
2\leq i \leq n-1\}$ and $\widetilde{A}=\{
log(x_{i_{1}}x_{i_{2}}...x_{i_{n}}) \ | \ i_{1}\in A_{1} ,i_{2}\in
A_{2},\ldots, i_{n-1}\in A_{n-1}, i_{n}\in A_{n}\}.$ Then the base
ring associated to transversal polymatroid presented by
$bf{\mathcal{\widetilde{A}}},$ \
$K[{\bf{\mathcal{\widetilde{A}}}}]$ is Gorenstein ring if and only
if there exists $\widetilde{N}\subset N$ such that
\[{\bf{R_{+}}}C_{\widetilde{A}}= \bigcap_{a\in \widetilde{N}}H^{+}_{a}\]
and $H^{+}_{a}$ with $a\in \widetilde{N}$ are just the facets of
the cone ${\bf{R_{+}}}C_{\widetilde{A}}.$

\end{conjecture}

\section{The description of some transversal polymatroids with Gorenstein base ring in dimension
$3$ and $4$.}

${\bf{Dimension\ 3.}}$

We consider the collection of sets ${\bf{\mathcal{A}}}=\{ \{1,2\},\{2,3\},\{3,1\} \}$.
The base ring associated to transversal polymatroid presented by ${\bf{\mathcal{A}}}$ is
\[ R = K[{\bf{\mathcal{A}}}]=K[x_1^{2}x_2, x_2^{2}x_1,  x_2^{2}x_3, x_3^{2}x_2, x_1^{2}x_3,
x_3^{2}x_1, x_1x_2x_3].\]

From { \cite{HH}}, $R$ is normal ring.

We can see $R = K[Q]$, where

$Q = \mathbb N\{(2,1,0),(1,2,0), (0,2,1),(0,1,2),(1,0,2),(2,0,1), (1,1,1)\}$.

Our aim is to describe the facets of $C=\mathbb R_{+}Q$.

It is easy to see that $C$ has $6$ facets, with the support planes given by the equations:
\[H_{\{1\}}: -x_{1}+2x_{2}+2x_{3}=0,\] \[ H_{\{2\}}: 2x_{1}-x_{2}+2x_{3}=0,\ \ \]
\[H_{\{3\}}:2x_{1}+2x_{2}-x_{3}=0,\ \ \]

\[H_{\{1,2\}}: x_{3}=0,\] \[H_{\{2,3\}}: x_{1}=0,\] \[H_{\{3,1\}}: x_{2}=0.\]
In fact, $C=H_{\{1\}}^{+}\cap H_{\{2\}}^{+}\cap H_{\{3\}}^{+}\cap
H_{\{1,2\}}^{+}\cap H_{\{2,3\}}^{+}\cap H_{\{3,1\}}^{+}.$ Since
$(1,1,1)$ is the only point in $ri(\mathbb R_{+}Q)$ then by
Danilov-Stanley theorem $R$ is a Gorenstein ring and
$\omega_{R}=R(-(1,1,1))$.

In order to obtain all the Gorenstein polymatroids of dimension
$3$, we will remove sequentially some facets of $C$. For instance,
if we remove the facet supported by $ H_{\{2\}}$ we obtain a new
cone $C'= H_{\{1\}}^{+}\cap H_{\{3\}}^{+}\cap H_{\{1,2\}}^{+}\cap
H_{\{2,3\}}^{+} \cap H_{\{3,1\}}^{+}.$ It is easy to note that
$C'=\mathbb R_{+}Q'$, where $Q'=Q + \mathbb N\{(0,3,0)\}$. $Q'$ is
a saturated semigroup, and moreover,
$K[Q']=K[{\bf{\mathcal{A}'}}]$, where ${\bf{\mathcal{A}'}} = \{
\{1,2\}, \{1,2,3\}, \{2,3\} \}$. The Danilov-Stanley theorem
assures us that $R'=K[Q']=K[{\bf{\mathcal{A}'}}]$ is still
Gorenstein with $\omega_{R'} = R'(-(1,1,1))$. (Remark. If we
remove the facet supported by $H_{\{3\}}$ or $H_{\{1\}}$, instead
of the facet supported by $H_{\{2\}}$ we obtain a new set
$\mathcal{A}'$ which is only a permutation of $1,2,3$.)

Suppose that we remove from $C'$ the facet supported by
$H_{\{3\}}$. We obtain a new cone $C''= H_{\{1\}}^{+}\cap
H_{\{1,2\}}^{+}\cap H_{\{2,3\}}^{+} \cap H_{\{3,1\}}^{+}.$ It is
easy to see that $C''= \mathbb R_{+}Q''$, where $Q''=Q' + \mathbb
N\{(0,0,3)\}$. $Q''$ is a saturated semigroup, and moreover,
$K[Q'']= K[{\bf{\mathcal{A}''}}]$, where ${\bf{\mathcal{A}''}} =
\{\{1,2,3\}, \{1,2,3\}, \{2,3\}\}$. The Danilov-Stanley theorem
implies that $R''=K[Q'']=K[{\bf{\mathcal{A}''}}]$ is Gorenstein
and $\omega_{R''} = R''(-(1,1,1))$. Finally, we remove from $C''$
the facet supported by $H_{\{1\}}$. We obtain the cone
$C'''=H_{\{1,2\}}^{+}\cap H_{\{2,3\}}^{+}\cap H_{\{3,1\}}^{+}$
which is a cone over $Q'''=Q'' + \mathbb N\{(3,0,0)\}$. $Q'''$ is
the saturated semigroup associated to the ring
$R'''=K[{\bf{\mathcal{A}'''}}]$, where
${\bf{\mathcal{A}'''}}=\{\{1,2,3\}, \{1,2,3\}, \{1,2,3\}\}$. Also,
$\omega_{R'''}=R'''(-(1,1,1))$.

Thus the base ring associated to the transversal polymatroids
presented by ${\bf{\mathcal{A}}}$, ${\bf{\mathcal{A}'}}$,
${\bf{\mathcal{A}''}}$, ${\bf{\mathcal{A}'''}}$ are Gorenstein
rings and for ${\bf{\mathcal{A}_{1}}}=\{ \{1,2\},\{2,3\} \}$ the
base ring presented by ${\bf{\mathcal{A}_{1}}}$ is the Segre
product $k[t_{11},t_{12}]*k[t_{21},t_{22}]$, thus is a Gorenstein
ring. All of them have dimension 3.\\

The computations made so far make us believe that all polymatroids
with Gorenstein base ring in dimension 3 are the ones classified
above.
\\
\[\]
${\bf{Dimension\ 4.}}$

We consider the collection of sets ${\bf{\mathcal{A}}}=\{ \{1,2\},\{2,3\},\{3,4\},\{4,1\} \}$.
The base ring associated to transversal polymatroid presented by ${\bf{\mathcal{A}}}$ is \[ R = K[{\bf{\mathcal{A}}}]=K[x_{i_{1}}x_{i_{2}}x_{i_{3}}x_{i_{4}}|\ i_{1}\in \{1,2\},i_{2}\in \{2,3\},i_{3}\in \{3,4\},i_{4}\in \{4,1\}] .\]

From { \cite{HH}}, $R$ is normal ring.

We can see $R = K[Q]$, where

\[Q = \mathbb N\{log(x_{i_{1}}x_{i_{2}}x_{i_{3}}x_{i_{4}})|\ i_{1}\in \{1,2\},i_{2}\in
\{2,3\},i_{3}\in \{3,4\},i_{4}\in \{4,1\} \}.\]

Our aim is to describe the facets of $C=\mathbb R_{+}Q$.
Using $\it Normaliz$ we have 12 facets of the cone $C=\mathbb R_{+}Q$:
\[ H_{\{1\}} : -x_{1}+x_{2}+x_{3}+x_{4}=0,\]
\[ H_{\{2\}} : x_{1}-x_{2}+x_{3}+x_{4}=0,\ \ \]
\[ H_{\{3\}} : x_{1}+x_{2}-x_{3}+x_{4}=0,\ \ \]
\[H_{\{4\}} : x_{1}+x_{2}+x_{3}-x_{4}=0,\ \ \]

 \[H_{\{1,2\}} : -x_{1}-x_{2}+3x_{3}+3x_{4}=0,\]
 \[ H_{\{2,3\}} : 3x_{1}-x_{2}-x_{3}+3x_{4}=0,\ \ \]
 \[ H_{\{3,4\}} : 3x_{1}+3x_{2}-x_{3}-x_{4}=0,\ \ \]
 \[ H_{\{1,4\}} : -x_{1}+3x_{2}+3x_{3}-x_{4}=0,\]

  \[H_{\{1,2,3\}} : x_{4}=0,\] \[ H_{\{2,3,4\}} : x_{1}=0,\] \[ H_{\{1,3,4\}} : x_{2}=0, \]
  \[H_{\{1,2,4\}} : x_{3}=0.\]

It is easy to see that $C=H_{\{1\}}^{+}\cap H_{\{2\}}^{+}\cap H_{\{3\}}^{+}\cap H_{\{4\}}^{+}\cap H_{\{1,2\}}^{+}\cap
H_{\{2,3\}}^{+}\cap H_{\{3,4\}}^{+} \cap H_{\{1,4\}}^{+}\cap H_{\{1,2,3\}}^{+}\cap H_{\{2,3,4\}}^{+}\cap H_{\{1,3,4\}}^{+}\cap H_{\{1,2,4\}}^{+}$. Since $(1,1,1,1)$ is the only point in $ri(\mathbb R_{+}Q)$ then by Danilov-Stanley theorem  $R$ is Gorenstein ring and $\omega_{R}=R(-(1,1,1,1))$.

Now I want to proceed like in dimension $3$ to give a large class
of transversal polymatroids with Gorenstein base ring. Using $\it
Normaliz$ we can give a complete description  modulo a permutation
of tha transversal polymatroids with Gorenstein base ring when we
start with ${\bf\mathcal{A}}=\{ \{1,2\},\{2,3\},\{3,4\},\{4,1\}
\}.$

For ${\bf{\mathcal{A}_{1}}}=\{ \{1,2,3\},\{2,3\},\{3,4\},\{4,1\} \}$ the associated cone is :\\
$C_{1}=H_{\{1\}}^{+}\cap H_{\{2\}}^{+}\cap H_{\{4\}}^{+}\cap H_{\{1,2\}}^{+}\cap
H_{\{2,3\}}^{+}\cap H_{\{1,4\}}^{+}\cap H_{\{1,2,3\}}^{+}\cap H_{\{2,3,4\}}^{+}\cap H_{\{1,3,4\}}^{+}\cap H_{\{1,2,4\}}^{+}.$\\

For ${\bf{\mathcal{A}_{2}}}=\{ \{1,2,3,4\},\{2,3\},\{3,4\},\{4,1\} \}$ the associated cone is :\\
$C_{2}=H_{\{1\}}^{+}\cap H_{\{2\}}^{+}\cap H_{\{1,2\}}^{+}\cap
H_{\{2,3\}}^{+}\cap H_{\{1,4\}}^{+}\cap H_{\{1,2,3\}}^{+}\cap H_{\{2,3,4\}}^{+}\cap H_{\{1,3,4\}}^{+}\cap H_{\{1,2,4\}}^{+}.$\\

For ${\bf{\mathcal{A}_{3}}}=\{ \{1,2,3,4\},\{2,3,4\},\{3,4\},\{4,1\} \}$ the associated cone is :\\
$C_{3}=H_{\{1\}}^{+}\cap H_{\{2\}}^{+}\cap H_{\{1,2\}}^{+}\cap
H_{\{2,3\}}^{+}\cap H_{\{1,2,3\}}^{+}\cap H_{\{2,3,4\}}^{+}\cap H_{\{1,3,4\}}^{+}\cap H_{\{1,2,4\}}^{+}.$\\

For ${\bf{\mathcal{A}_{4}}}=\{ \{1,2,3,4\},\{1,2,3,4\},\{3,4\},\{4,1\} \}$ the associated cone is :\\
$C_{4}=H_{\{2\}}^{+}\cap H_{\{1,2\}}^{+}\cap
H_{\{2,3\}}^{+}\cap H_{\{1,2,3\}}^{+}\cap H_{\{2,3,4\}}\cap H_{\{1,3,4\}}\cap H_{\{1,2,4\}}.$\\

For ${\bf{\mathcal{A}_{5}}}=\{ \{1,2,3,4\},\{1,2,3,4\},\{1,3,4\},\{4,1\} \}$ the associated cone is :\\
$C_{5}=H_{\{2\}}^{+}\cap
H_{\{2,3\}}^{+}\cap H_{\{1,2,3\}}^{+}\cap H_{\{2,3,4\}}^{+}\cap H_{\{1,3,4\}}^{+}\cap H_{\{1,2,4\}}^{+}.$\\

For ${\bf{\mathcal{A}_{6}}}=\{ \{1,2,3,4\},\{1,2,3,4\},\{1,2,3,4\},\{4,1\} \}$ the associated cone is :\\
$C_{6}=H_{\{2,3\}}^{+}\cap H_{\{1,2,3\}}^{+}\cap H_{\{2,3,4\}}^{+}\cap H_{\{1,3,4\}}^{+}\cap H_{\{1,2,4\}}^{+}.$\\

For ${\bf{\mathcal{A}_{7}}}=\{ \{1,2,3,4\},\{1,2,3,4\},\{1,2,3,4\},\{1,2,3,4\} \}$ the associated cone is :\\
$C_{7}=H_{\{1,2,3\}}^{+}\cap H_{\{2,3,4\}}^{+}\cap H_{\{1,3,4\}}^{+}\cap H_{\{1,2,4\}}^{+}.$\\

For ${\bf{\mathcal{A}_{8}}}=\{ \{1,2,3\},\{1,2,3\},\{3,4\},\{4,1\} \}$ the associated cone is :\\
$C_{8}=H_{\{2\}}^{+}\cap H_{\{4\}}^{+}\cap H_{\{1,2\}}^{+}\cap
H_{\{2,3\}}^{+}\cap H_{\{1,2,3\}}^{+}\cap H_{\{2,3,4\}}^{+}\cap H_{\{1,3,4\}}^{+}\cap H_{\{1,2,4\}}^{+}.$\\

For ${\bf{\mathcal{A}_{9}}}=\{ \{1,2,3\},\{1,2,3\},\{1,3,4\},\{4,1\} \}$ the associated cone is :\\
$C_{9}=H_{\{2\}}^{+}\cap H_{\{4\}}^{+}\cap
H_{\{2,3\}}^{+}\cap H_{\{1,2,3\}}^{+}\cap H_{\{2,3,4\}}^{+}\cap H_{\{1,3,4\}}^{+}\cap H_{\{1,2,4\}}^{+}.$\\

For ${\bf{\mathcal{A}_{10}}}=\{ \{1,2,3\},\{1,2,3\},\{1,2,3,4\},\{4,1\} \}$ the associated cone is :\\
$C_{10}=H_{\{4\}}^{+}\cap
H_{\{2,3\}}^{+}\cap H_{\{1,2,3\}}^{+}\cap H_{\{2,3,4\}}^{+}\cap H_{\{1,3,4\}}^{+}\cap H_{\{1,2,4\}}^{+}.$\\

For ${\bf{\mathcal{A}_{11}}}=\{ \{1,2,3\},\{1,2,3\},\{1,3,4\},\{1,3,4\} \}$ the associated cone is :\\
$C_{11}=H_{\{2\}}^{+}\cap H_{\{4\}}^{+}\cap H_{\{1,2,3\}}^{+}\cap H_{\{2,3,4\}}^{+}\cap H_{\{1,3,4\}}^{+}\cap H_{\{1,2,4\}}^{+}.$\\

For ${\bf{\mathcal{A}_{12}}}=\{ \{1,2,3\},\{1,2,3\},\{1,3,4\},\{1,3,4\} \}$ the associated cone is :\\
$C_{12}=H_{\{4\}}^{+}\cap H_{\{1,2,3\}}^{+}\cap
H_{\{2,3,4\}}^{+}\cap H_{\{1,3,4\}}^{+}\cap H_{\{1,2,4\}}^{+}.$\\

The next four examples of transversal polymatroids with Gorenstein
base ring are different like above.

For ${\bf{\mathcal{A}_{13}}}=\{ \{1,2,3\},\{2,3,4\}\}$ the Hilbert
series of base ring $K[{\bf{\mathcal{A}_{13}}}]$ is:\\
$H_{K[{\bf{\mathcal{A}_{13}}}]}(t)=\frac{1+4t+t^{2}}{(1-t)^{4}}.$

For ${\bf{\mathcal{A}_{14}}}=\{ \{1,2,3,4\},\{2,3,4\}\}$ the
Hilbert series of base ring $K[{\bf{\mathcal{A}_{14}}}]$ is:\\
$H_{K[{\bf{\mathcal{A}_{14}}}]}(t)=\frac{1+5t+t^{2}}{(1-t)^{4}}.$

For ${\bf{\mathcal{A}_{15}}}=\{ \{1,2,3,4\},\{1,2,3,4\}\}$ the
Hilbert series of base ring $K[{\bf{\mathcal{A}_{15}}}]$ is:\\
$H_{K[{\bf{\mathcal{A}_{15}}}]}(t)=\frac{1+6t+t^{2}}{(1-t)^{4}}.$

For ${\bf{\mathcal{A}_{16}}}=\{ \{1,2\},\{2,3\},\{3,4\}\}$ the
Hilbert series of base ring $K[{\bf{\mathcal{A}_{16}}}]$ is:\\
$H_{K[{\bf{\mathcal{A}_{16}}}]}(t)=\frac{1+4t+t^{2}}{(1-t)^{4}}.$

It seems that also in dimension 4 our examples cover all
transversal polymatroids with Gorenstein base ring.

{\em \noindent University of Ploie\c sti, Department of Mathematics\\[0pt]
Bd. Bucure\c sti, 39, Ploie\c sti, Romania}\\[0pt]
e-mail: nastefan@mail.upg-ploiesti.ro

\end{document}